\documentclass{amsart}
\usepackage{wasowicz}

\newtheorem*{ThA}{Theorem A}

\info{Published: Publ. Math. Debrecen 68 (2006), 171--182}

\begin{document}
\author{\SW}
\address{\SWaddr}
\email{\SWmail}
\title{Some properties of generalized higher-order convexity}
\date{November 18, 2004}
\begin{abstract}
 The generalized divided differences are introduced. They are applied
 to investigate some properties characterizing generalized
 higher-order convexity. Among others some support-type property
 is proved.
\end{abstract}
\maketitle
\section{Introduction}
Let $I\subset\R$ be an interval and let
$\omega_1,\dots,\omega_n:I\to\R$ be continuous functions. For $n$
distinct points $x_{i_1},\dots,x_{i_n}\in I$ we define
\begin{equation}\label{Vn}
 V_n(x_{i_1},\dots,x_{i_n})
 =\begin{vmatrix}
   \omega_1(x_{i_1})&\dots&\omega_1(x_{i_n})\\
   \vdots&\ddots&\vdots\\
   \omega_n(x_{i_1})&\dots&\omega_n(x_{i_n})
  \end{vmatrix}.
\end{equation}
A system $\pmb{\omega}=(\omega_1,\dots,\omega_n)$ is called
a~\emph{Chebyshev system on~$I$} if $V_n(x_1,\dots,$ $x_n)\ne 0$
for any $x_1,\dots,x_n\in I$ such that $x_1<\dots<x_n$.
\begin{exmp}\label{Ex_Ch}
 The systems $\pmb{\omega}=(1,x,\dots,x^{n-1})$,
 $\pmb{\omega}=(\e^{\alpha_1x},\dots,\e^{\alpha_nx})$
 (for any distinct $\alpha_1,\dots,\alpha_n\in\R$) are
 Chebyshev systems on any interval.
\end{exmp}
\begin{rem}\label{n-parameter}
By the Cramer Rule a~linear span of a~Chebyshev system
$\pmb{\omega}=(\omega_1,\dots,\omega_n)$
is an $n$-parameter family on~$I$, i.e. for any $n$ distinct points
$x_1,\dots,x_n\in I$ and for any $y_1,\dots,y_n\in\R$ there exists
exactly one function $\omega=c_1\omega_1+\dots+c_n\omega_n$
(where $c_1,\dots,c_n\in\R$ are the constants) such that
$\omega(x_i)=y_i$, $i=1,\dots,n$. Such families were considered by
Tornheim~\cite{Tor50} (see also Beckenbach~\cite{Bec37},
Beckenbach and Bing~\cite{BecBin45}).
\end{rem}
If $\pmb{\omega}=(\omega_1,\dots,\omega_n)$ is a~Chebyshev system
on~$I$ then by continuity of $\omega_1,\dots,\omega_n$ the determinant
$V_n(x_1,\dots,x_n)$ does not change the sign in a~connected set
$\bigl\{(x_1,\dots,x_n)\in I:x_1<\dots<x_n\bigr\}$. Then a~Chebyshev
system~$\pmb{\omega}$ is called \emph{positive} (\emph{negative}) if
$V_n(x_1,\dots,x_n)>0$ ($V_n(x_1,\dots,x_n)<0$) for all
$x_1,\dots,x_n\in I$ such that $x_1<\dots<x_n$. Notice that the
Chebyshev systems of Example~\ref{Ex_Ch} are positive.
\begin{rem}
 Throughout the paper we will often assume that
 $\pmb{\omega}=(\omega_1,\dots,\omega_n)$ is such a~Chebyshev system
 on~$I$ that $(\omega_1,\dots,\omega_{n-1})$ is also a~Chebyshev
 system on~$I$. This assumption is not too restrictive. Many Chebyshev
 systems have this property, e.g. the systems mentioned in
 Example~\ref{Ex_Ch}. However $(\cos x,\sin x)$ is a~Chebyshev system
 on $(0,\pi)$ but $(\cos x)$ is not a~Chebyshev system on $(0,\pi)$.
 \par
 We will also assume that $\pmb{\omega}=(\omega_1,\dots,\omega_n)$ is
 a~positive Chebyshev system on~$I$ such that
 $(\omega_1,\dots,\omega_{n-1})$ is also a~positive Chebyshev system
 on~$I$. The systems of Example~\ref{Ex_Ch} satisfy this assumption
 as well. But there are Chebyshev systems which do not have this
 property. Notice that $(-1,-x)$ is a~positive Chebyshev system on
 any interval but $(-1)$ is a~negative one.
\end{rem}
For a~function $f:I\to\R$ and for $n+1$ distinct points
$x_1,\dots,x_{n+1}\in I$ we define
\begin{equation}\label{Dn}
 D_n(x_1,\dots,x_{n+1};f)
 =\begin{vmatrix}
   \omega_1(x_1)&\dots&\omega_1(x_{n+1})\\
   \vdots&\ddots&\vdots\\
   \omega_n(x_1)&\dots&\omega_n(x_{n+1})\\
   f(x_1)&\dots&f(x_{n+1})
  \end{vmatrix}.
\end{equation}
Let $\pmb{\omega}=(\omega_1,\dots,\omega_n)$ be a~Chebyshev
system on~$I$. A~function $f:I\to\R$ is called
$\pmb{\omega}$-$n$-\emph{convex} if for any $n$ distinct points
$x_1,\dots,x_n\in I$ such that $x_1<\dots<x_n$ the (uniquely
determined) function $\omega=c_1\omega_1+\dots+c_n\omega_n$ such that
$\omega(x_i)=f(x_i)$, $i=1,\dots,n$, fulfils the conditions
\begin{align*}
 (-1)^n\bigl(f(x)-\omega(x)\bigr)&\ge 0\quad\text{for }x\le x_1,\\
 (-1)^{n+i}\bigl(f(x)-\omega(x)\bigr)&\ge 0
 \quad\text{for }x_i\le x\le x_{i+1},\;i=1,\dots,n-1,\\
 f(x)-\omega(x)&\ge 0\quad\text{for }x\ge x_n
\end{align*}
(see~\cite{BesPal04}, \cite{Tor50}; for $\pmb{\omega}$-$n$-convexity
with respect to $\pmb{\omega}=(1,x,\dots,x^{n-1})$ see also
\cite{Kuc85}, \cite{RobVar73}).
\par
Observe that for $n=2$ and $\pmb{\omega}=(1,x)$
$\pmb{\omega}$-$2$-convexity reduces to convexity in the usual
sense. Indeed, $f$ is $\pmb{\omega}$-$2$-convex if and only if for
any $x_1,x_2\in I$ such that $x_1<x_2$ there exists an affine function
$\omega(x)=c_1+c_2x$, $x\in I$, such that $\omega(x_i)=f(x_i)$,
$i=1,2$ and $f\le\omega$ on $[x_1,x_2]$ (and $\omega\le f$ on
$I\setminus[x_1,x_2]$). This statement is evidently equivalent to
convexity of~$f$.

For $\pmb{\omega}=(1,x,\dots,x^{n-1})$ $\pmb{\omega}$-$n$-convex
functions are convex functions of higher orders (see~\cite{Kuc85},
\cite{Pop44}, \cite{Pop34}, \cite{RobVar73}, \cite{Tor50}).
\par
Bessenyei and P\'ales obtained the following result
(\cite[Theorem 2 $(i)\Leftrightarrow(iii)$]{BesPal04}).
\begin{ThA}
 Let $\pmb{\omega}=(\omega_1,\dots,\omega_n)$ be a~positive Chebyshev
 system on~$I$. A function~$f:I\to\R$ is $\pmb{\omega}$-$n$-convex if
 and only if
 \[
  D_n(x_1,\dots,x_{n+1};f)\ge 0
 \] for all
 $x_1,\dots,x_{n+1}\in I$ such that $x_1<\dots<x_{n+1}$.
\end{ThA}
N\"orlund~\cite{Nor26} considered the divided differences given by
the following recurrence: 
\begin{equation}\label{div_diff_classic}
 [x_1,f]=f(x_1)\quad\text{and}\quad
 [x_1,\dots,x_{n+1};f]
 =\frac{[x_2,\dots,x_{n+1};f]-[x_1,\dots,x_n;f]}{x_{n+1}-x_1}
\end{equation}
(cf. also~\cite{Kuc85}, \cite{Pop34}, \cite{RobVar73}).
Now we are going to generalize this notion.
\par
Let $\pmb{\omega}=(\omega_1,\dots,\omega_n)$ be a~Chebyshev system
on~$I$ such that $(\omega_1,\dots,$ $\omega_{n-1})$ is also a~Chebyshev
system on~$I$. For $n$ distinct points $x_1,\dots,x_n\in I$ we
introduce the \emph{generalized divided differences} by the formula
\begin{equation}\label{div_diff}
 [x_1,\dots,x_n;f]_{\pmb{\omega}}
 =\frac{D_{n-1}(x_1,\dots,x_n;f)}{V_n(x_1,\dots,x_n)}.
\end{equation}
For $\pmb{\omega}=(1,x,\dots,x^{n-1})$ the generalized divided
difference $[x_1,\dots,x_n;f]_{\pmb{\omega}}$ is equal to
$[x_1,\dots,x_n;f]$ given by~\eqref{div_diff_classic}
(see \cite{Kuc85}, \cite{Pop34}).

\begin{rem}\label{Rem1}
 The generalized divided differences are symmetric. Na\-me\-ly, if
 $(x_{i_1},\dots,x_{i_n})$ is a~permutation of $(x_1,\dots,x_n)$ then
 \begin{equation}\label{perm}
  [x_1,\dots,x_n;f]_{\pmb{\omega}}
  =[x_{i_1},\dots,x_{i_n};f]_{\pmb{\omega}}.
 \end{equation}
 This is a~simple consequence of the properties of determinants. To
 get $[x_{i_1},\dots,x_{i_n};f]_{\pmb{\omega}}$ we need to make the
 same inversions both in the numerator and in the denominator of
 $[x_1,\dots,x_n;f]_{\pmb{\omega}}$.
\end{rem}
In this paper we prove in Theorem~\ref{Th1} an analogue
of~\eqref{div_diff_classic} for generalized divided differences,
which seems to be very convenient to investigate the properties of
$\pmb{\omega}$-$n$-convexity.
Using Theorem~\ref{Th1} we prove in Theorem~\ref{Th2} that
a~function~$f$ is $\pmb{\omega}$-$n$-convex if and only if its
generalized divided differences are nondecreasing. Another
characterization of $\pmb{\omega}$-$n$-convexity is some support-type
property proved in Theorem~\ref{Th3}. The classical support theorems
state that for a~real function~$f$ and for some element $x_0$ of its
domain under suitable assumptions there exists a~function $g$ (the
supporting function) such that $g(x_0)=f(x_0)$ and $g\le f$. Our
Theorem~\ref{Th3} is not the classical support theorem. The graph of
obtained ''supporting function'' meets the graph of the ''supported
function''~$f$ at~$n-1$ points $x_1<\dots<x_{n-1}$ and passing through
$x_1,\dots,x_{n-2}$ it changes successively the side of the graph
of~$f$ being the classical supporting function in the subinterval
$(x_{n-2},+\infty)\cap I$. It is worth mentioning that this result
extends the recent result of Bessenyei and
P\'ales~(\cite[Theorem~4 $(i)\Leftrightarrow(iii)$]{BesPal03}) concerning
$\pmb{\omega}$-2-convexity.

\section{Some property of generalized divided differences}

We start with the generalization of~\eqref{div_diff_classic}. This is
an equation~\eqref{dd} below which seems to be very convenient to
investigate the properties of $\pmb{\omega}$-$n$-convexity. It is easy
to observe that for $\pmb{\omega}=(1,x,\dots,x^{n-1})$~\eqref{dd}
reduces to~\eqref{div_diff_classic}.

\begin{thm}\label{Th1}
 Let $n\ge 2$, let $\pmb{\omega}=(\omega_1,\dots,\omega_n)$ be
 a~Chebyshev system on~$I$ such that $(\omega_1,\dots,\omega_{n-1})$
 is also a~Chebyshev system on~$I$ and let $f:I\to\R$. Then
 \begin{equation}\label{dd}
  [x_2,\dots,x_{n+1};f]_{\pmb{\omega}}
  -[x_1,\dots,x_n;f]_{\pmb{\omega}}
  =\frac{D_n(x_1,\dots,x_{n+1};f)V_{n-1}(x_2,\dots,x_n)}
        {V_n(x_2,\dots,x_{n+1})V_n(x_1,\dots,x_n)}
 \end{equation}
 for any $n+1$ distinct points $x_1,\dots,x_{n+1}\in I$.
\end{thm}
\begin{proof}
Since $(\omega_1,\dots,\omega_{n-1})$ is a Chebyshev system then by
Remark~\ref{n-parameter} we can choose the constants
$c_1,\dots,c_{n-1}$ such that for
$\omega=c_1\omega_1+\dots+c_{n-1}\omega_{n-1}$ we have
$\omega(x_k)=f(x_k)$, $k=2,\dots,n$. Then for $f^{\ast}=f-\omega$
we obtain
\begin{equation}\label{Th1_1}
 f^{\ast}(x_2)=\dots=f^{\ast}(x_n)=0.
\end{equation}
By the elementary properties of determinants we get
$[x_2,\dots,x_{n+1};\omega]_{\pmb{\omega}}=0$ and
\[
 [x_2,\dots,x_{n+1};f]_{\pmb{\omega}}=
 [x_2,\dots,x_{n+1};\omega+f^{\ast}]_{\pmb{\omega}}=
 [x_2,\dots,x_{n+1};f^{\ast}]_{\pmb{\omega}}.
\]
Similarly
\begin{align*}
 [x_1,\dots,x_n;f]_{\pmb{\omega}}&=[x_1,\dots,x_n;f^{\ast}]_{\pmb{\omega}}\\
 \intertext{and}
 D_n(x_1,\dots,x_{n+1};f)&=D_n(x_1,\dots,x_{n+1};f^{\ast}).
\end{align*}
Then replacing in~\eqref{dd} $f$ by $\omega+f^{\ast}$ and using the
previous three equations we can see that it
is enough to prove~\eqref{dd} only for~$f^{\ast}$.
\par
Expanding $D_n(x_1,\dots,x_{n+1};f^{\ast})$ by its last row and
using~\eqref{Th1_1} we obtain
\begin{multline}\label{Th1_2}
 D_n(x_1,\dots,x_{n+1};f^{\ast})\\
 =(-1)^nf^{\ast}(x_1)V_n(x_2,\dots,x_{n+1})+
  f^{\ast}(x_{n+1})V_n(x_1,\dots,x_n).
\end{multline}
By~\eqref{div_diff} we have
\begin{multline*}
 [x_2,\dots,x_{n+1};f^{\ast}]_{\pmb{\omega}}
  -[x_1,\dots,x_n;f^{\ast}]_{\pmb{\omega}}\\
 =\frac{D_{n-1}(x_2,\dots,x_{n+1};f^{\ast})}{V_n(x_2,\dots,x_{n+1})}-
  \frac{D_{n-1}(x_1,\dots,x_n;f^{\ast})}{V_n(x_1,\dots,x_n)}.
\end{multline*}
Expanding the numerators by the last rows and using~\eqref{Th1_1}
we get
\begin{align*}
 &[x_2,\dots,x_{n+1};f^{\ast}]_{\pmb{\omega}}
  -[x_1,\dots,x_n;f^{\ast}]_{\pmb{\omega}}\\
 &=\frac{f^{\ast}(x_{n+1})V_{n-1}(x_2,\dots,x_n)}{V_n(x_2,\dots,x_{n+1})}
 -\frac{(-1)^{n+1}f^{\ast}(x_1)V_{n-1}(x_2,\dots,x_n)}{V_n(x_1,\dots,x_n)}
\end{align*}
Then by~\eqref{Th1_2} we obtain~\eqref{dd} for $f^{\ast}$ which
finishes the proof.
\end{proof}

\section{Some characterizations of $\pmb{\omega}$-$n$-convexity}

\begin{cor}\label{Co1}
 Let $n\ge 2$, let $\pmb{\omega}=(\omega_1,\dots,\omega_n)$ be
 a~positive Chebyshev system on~$I$ such that
 $(\omega_1,\dots,\omega_{n-1})$ is also a~positive Chebyshev system
 on~$I$. A~function $f:I\to\R$ is $\pmb{\omega}$-$n$-convex if and
 only if
 \[
  [x_2,\dots,x_{n+1};f]_{\pmb{\omega}}\ge
  [x_1,\dots,x_n;f]_{\pmb{\omega}}
 \]
 for all $x_1,\dots,x_{n+1}\in I$ such that $x_1<\dots<x_{n+1}$.
\end{cor}
\begin{proof}
 Since $\pmb{\omega}$ and $(\omega_1,\dots,\omega_{n-1})$ are positive
 Chebyshev systems then the determinants $V_n(x_1,\dots,x_n)$,
 $V_n(x_2,\dots,x_{n+1})$ and $V_{n-1}(x_2,\dots,$ $x_n)$ are positive
 for all $x_1,\dots,$ $x_{n+1}\in I$ such that $x_1<\dots<x_{n+1}$.
 Then Corollary~\ref{Co1} follows immediately by~\eqref{dd} and by
 Theorem~$A$.
\end{proof}
\begin{rem}
 Corollary~\ref{Co1} generalizes the equivalence
 $(i)\Leftrightarrow(ii)$ of Theorem~4 of~\cite{BesPal03}.
 We obtain it using Corollary~\ref{Co1} for $n=2$.
\end{rem}

Next we state that a~function~$f$ is $\pmb{\omega}$-$n$-convex if
and only if its generalized divided differences are nondecreasing.
For $n=2$ and $\pmb{\omega}=(1,x)$ we obtain the very well known
characterization of the usual convexity: a~function~$f$ is convex if
and only if its difference quotients are nondecreasing. By $I^0$ we
denote the interior of~$I$.

\begin{thm}\label{Th2}
 Let $n\ge 2$, let $\pmb{\omega}=(\omega_1,\dots,\omega_n)$ be
 a~positive Chebyshev system on~$I$ such that
 $(\omega_1,\dots,\omega_{n-1})$ is also a~positive Chebyshev system
 on~$I$. A~function $f:I\to\R$ is $\pmb{\omega}$-$n$-convex if and
 only if for all $x_1,\dots,x_{n-1}\in I^0$ such that
 $x_1<\dots<x_{n-1}$ the~function
 $x\mapsto[x_1,\dots,x_{n-1},x;f]_{\pmb{\omega}}$ is nondecreasing
 on the set $I\setminus\{x_1,\dots,x_{n-1}\}$.
\end{thm}

\begin{proof}
 Take $x_1,\dots,x_{n-1}\in I^0$ such that $x_1<\dots<x_{n-1}$
 and $x,y\in I\setminus\{x_1,\dots,x_{n-1}\}$ such that $x<y$.
 The points $x_1,\dots,x_{n-1}$ divide the set
 $I\setminus\{x_1,\dots,x_{n-1}\}$ into $n$ subintervals
 $I_1=(-\infty,x_1)\cap I$, $I_s=(x_{s-1},x_s)$, $s=2,\dots,n-1$
 (if $n\ge 3$) and $I_n=(x_{n-1},+\infty)\cap I$. Let $x\in I_j$,
 $y\in I_k$. Since $x<y$ then $j\le k$. There are $j-1$~inversions
 of~$x$ needed to transform the ordered system of~$n$ points
 $(x_1,\dots,x,\dots,x_{n-1})$ to the system $(x,x_1,\dots,x_{n-1})$.
 Then
 \begin{equation}\label{Th2_1}
  V_n(x,x_1,\dots,x_{n-1})=(-1)^{j-1}V_n(x_1,\dots,x,\dots,x_{n-1}).
 \end{equation}
 We need $n-k$ inversions of~$y$ to transform the ordered system
 of~$n$ points $(x_1,\dots,y,$ $\dots,x_{n-1})$ to the system
 $(x_1,\dots,x_{n-1},y)$. Then
 \begin{equation}\label{Th2_2}
  V_n(x_1,\dots,x_{n-1},y)=(-1)^{n-k}V_n(x_1,\dots,y,\dots,x_{n-1}).
 \end{equation}
 Observe that starting from the ordered system of~$n+1$ points
 $(x_1,\dots,x,\dots,$ $y,\dots,x_{n-1})$ after $j-1$~inversions of~$x$
 and $n-k$~inversions of~$y$ we get the system
 $(x,x_1,\dots,x_{n-1},y)$. Then
 \begin{multline}\label{Th2_3}
  D_n(x,x_1,\dots,x_{n-1},y;f)\\
  =(-1)^{j-1+n-k}D_n(x_1,\dots,x,\dots,y,\dots,x_{n-1};f).
 \end{multline}
 By \eqref{Th2_1},~\eqref{Th2_2},~\eqref{Th2_3}, Remark~\ref{Rem1} and
 Theorem~\ref{Th1} we obtain
 \begin{align*}
  [x_1,\dots,&x_{n-1},y;f]_{\pmb{\omega}}
  -[x_1,\dots,x_{n-1},x;f]_{\pmb{\omega}}\\
  &=[x_1,\dots,x_{n-1},y;f]_{\pmb{\omega}}
  -[x,x_1,\dots,x_{n-1};f]_{\pmb{\omega}}\\
  &=\frac{D_n(x,x_1,\dots,x_{n-1},y;f)V_{n-1}(x_1,\dots,x_{n-1})}
        {V_n(x_1,\dots,x_{n-1},y)V_n(x,x_1,\dots,x_{n-1})}\\
  &=\frac{D_n(x_1,\dots,x,\dots,y,\dots,x_{n-1})V_{n-1}(x_1,\dots,x_{n-1})}
        {V_n(x_1,\dots,y,\dots,x_{n-1})V_n(x_1,\dots,x,\dots,x_{n-1})}\,.
 \end{align*}
 Observe that the determinants $V_{n-1}(x_1,\dots,x_{n-1})$,
 $V_n(x_1,\dots,y,\dots,x_{n-1})$ and $V_n(x_1,$ $\dots,x,\dots,x_{n-1})$
 are positive since $\pmb{\omega}$ and $(\omega_1,\dots,\omega_{n-1})$
 are positive Chebyshev systems and the systems of points involved are
 ordered. Then Theorem~\ref{Th2} follows immediately by Theorem~A.
\end{proof}

\section{Support-type property of $\pmb{\omega}$-$n$-convexity}

In this section we are going to prove some kind of support theorem. In
the classical approach the graph of the supporting function lies below
(precisely not above) the graph of the supported function and it meets
this graph (at least) at one point. For a~discussion of our approach
see the Introduction. The ''support'' property proved in
Theorem~\ref{Th3} characterizes $\pmb{\omega}$-$n$-convexity. Let
us mention that Ger \cite[Corollary 2]{Ger94} proved the classical
support theorem for convex functions of an odd order~$n$. Here the
supporting function is the polynomial of an order at most~$n$. The
classical polynomial support property is no longer valid for the
convex functions of an even order (see \cite[Remark 1]{Ger94}). Our
Theorem~\ref{Th3} (applied for $\pmb{\omega}=(1,x,\dots,x^{n-1})$)
characterizes the convexity of both odd and even order. We start with
the following technical result.

\begin{lem}\label{Lm1}
 Let $n\ge 2$, let $\pmb{\omega}=(\omega_1,\dots,\omega_n)$ be
 a~Chebyshev system on~$I$ such that $(\omega_1,\dots,\omega_{n-1})$
 is also a~Chebyshev system on~$I$, let $c_n\in\R$ and let $f:I\to\R$.
 Then for any $n-1$ distinct points $x_1,\dots,x_{n-1}\in I^0$ there
 exist the constants $c_1,\dots,c_{n-1}\in\R$ such that for
 $\omega=c_1\omega_1+\dots+c_{n-1}\omega_{n-1}+c_n\omega_n$ we have
 $\omega(x_k)=f(x_k)$, $k=1,\dots,n-1$ and
 \[
  f(x)-\omega(x)=
  \frac{D_{n-1}(x_1,\dots,x_{n-1},x;f)-c_nV_n(x_1,\dots,x_{n-1},x)}
       {V_{n-1}(x_1,\dots,x_{n-1})}
 \]
 for all $x\in I\setminus\{x_1,\dots,x_{n-1}\}$.
\end{lem}
\begin{proof}
 Fix $c_n\in\R$. Since $(\omega_1,\dots,\omega_{n-1})$ is a~Chebyshev
 system, the constants $c_1,\dots,$ $c_{n-1}$ are (uniquely) determined
 by the system of linear equations
 \[
  c_1\omega_1(x_k)+\dots+c_{n-1}\omega_{n-1}(x_k)
  =f(x_k)-c_n\omega_n(x_k),\quad k=1,\dots,n-1.
 \]
 Then for $\omega=c_1\omega_1+\dots+c_{n-1}\omega_{n-1}+c_n\omega_n$
 we have
 \begin{equation}\label{Lm1_1}
  \omega(x_k)=f(x_k),\quad k=1,\dots,n-1.
 \end{equation}
 Let $x\in I\setminus\{x_1,\dots,x_{n-1}\}$.
 Expanding the determinant $D_{n-1}(x_1,\dots,x_{n-1},$ $x;f-\omega)$
 by the last row and using~\eqref{Lm1_1} we get
 \begin{equation}\label{Lm1_2}
  D_{n-1}(x_1,\dots,x_{n-1},x;f-\omega)
  =\bigl(f(x)-\omega(x)\bigr)V_{n-1}(x_1,\dots,x_{n-1}).
 \end{equation}
 Since $D_{n-1}(x_1,\dots,x_{n-1},x;\omega_k)=0$, $k=1,\dots,n-1$,
 then
 \begin{multline}\label{Lm1_3}
  D_{n-1}(x_1,\dots,x_{n-1},x;\omega)
  =\sum_{k=1}^{n-1}c_kD_{n-1}(x_1,\dots,x_{n-1},x;\omega_k)\\
  +c_nD_{n-1}(x_1,\dots,x_{n-1},x;\omega_n)
  =c_nV_n(x_1,\dots,x_{n-1},x).
 \end{multline}
 Then using \eqref{Lm1_2} and~\eqref{Lm1_3} we obtain
 \begin{align*}
  &f(x)-\omega(x)=\frac{D_{n-1}(x_1,\dots,x_{n-1},x;f-\omega)}{V_{n-1}(x_1,\dots,x_{n-1})}\\
  &=\frac{D_{n-1}(x_1,\dots,x_{n-1},x;f)-D_{n-1}(x_1,\dots,x_{n-1},x;\omega)}{V_{n-1}(x_1,\dots,x_{n-1})}\\
  &=\frac{D_{n-1}(x_1,\dots,x_{n-1},x;f)-c_nV_n(x_1,\dots,x_{n-1},x)}{V_{n-1}(x_1,\dots,x_{n-1})},
 \end{align*}
 which was to be proved.
\end{proof}

Next we prove the support-type result mentioned at the beginning of
this section.

\begin{thm}\label{Th3}
 Let $n\ge 2$, let $\pmb{\omega}=(\omega_1,\dots,\omega_n)$ be
 a~positive Chebyshev system on~$I$ such that
 $(\omega_1,\dots,\omega_{n-1})$ is also a~positive Chebyshev system
 on~$I$. A~function $f:I\to\R$ is $\pmb{\omega}$-$n$-convex if and
 only if for all $x_1,\dots,x_{n-1}\in I^0$ such that
 $x_1<\dots<x_{n-1}$ there exist the constants $c_1,\dots,c_n\in\R$
 such that for $\omega=c_1\omega_1+\dots+c_n\omega_n$ we have
 $\omega(x_k)=f(x_k)$, $k=1,\dots,n-1$ and
 \begin{align}
  \label{Th3_i}
  (-1)^{n-1}\bigl(f(x)-\omega(x)\bigr)&\le 0\quad\text{for $x\in I$ such that $x<x_1$},\\
  \label{Th3_ii}
  (-1)^{n-k}\bigl(f(x)-\omega(x)\bigr)&\le 0\quad
  \text{for }x_{k-1}<x<x_k,\;k=2,\dots,n-1,\\
  \label{Th3_iii}
  f(x)-\omega(x)&\ge 0\quad\text{for $x\in I$ such that $x>x_{n-1}$}
 \end{align}
 (for $n=2$ there are no inequalities~\eqref{Th3_ii}).
\end{thm}
\begin{proof}
 Assume that $f$ is $\pmb{\omega}$-$n$-convex and fix
 $x_1,\dots,x_{n-1}\in I^0$ such that $x_1<\dots<x_{n-1}$. By
 Theorem~\ref{Th2} the function
 $x\mapsto[x_1,\dots,$ $x_{n-1},x;f]_{\pmb{\omega}}$
 is nondecreasing on the set $I\setminus\{x_1,\dots,x_{n-1}\}$.
 Then we define
 \begin{equation}\label{Th3_1}
  c_n=\lim\limits_{x\to x_{n-1}^+}
       [x_1,\dots,x_{n-1},x;f]_{\pmb{\omega}}.
 \end{equation}
 By Lemma~\ref{Lm1} there exist the constants
 $c_1,\dots,c_{n-1}\in\R$ such that for
 $\omega=c_1\omega_1+\dots+c_{n-1}\omega_{n-1}+c_n\omega_n$ we have
 $\omega(x_k)=f(x_k)$, $k=1,\dots,n-1$. Then to prove the necessity
 we have to check the inequalities~\eqref{Th3_i},~\eqref{Th3_ii}
 and \eqref{Th3_iii}. We start with~\eqref{Th3_iii}. Fix $x\in I$
 such that $x>x_{n-1}$.
 Theorem~\ref{Th2} and~\eqref{Th3_1} yield
 $c_n\le[x_1,\dots,x_{n-1},x;f]_{\pmb{\omega}}$. Then
 by~\eqref{div_diff} we have
 \[
  c_n\le\frac{D_{n-1}(x_1,\dots,x_{n-1},x;f)}{V_n(x_1,\dots,x_{n-1},x)}.
 \]
 Since $x_1<\dots<x_{n-1}<x$, then $V_n(x_1,\dots,x_{n-1},x)>0$,
 whence
 \[
  D_{n-1}(x_1,\dots,x_{n-1},x;f)-c_nV_n(x_1,\dots,x_{n-1},x)\ge 0.
 \]
 Dividing both sides of this inequality by
 $V_{n-1}(x_1,\dots,x_{n-1})>0$ and using Lemma~\ref{Lm1} we obtain
 $f(x)-\omega(x)\ge 0$.
 \par
 Let us now check~\eqref{Th3_i} and~\eqref{Th3_ii}. Similarly as in
 the proof of Theorem~\ref{Th2} denote $I_1=(-\infty,x_1)\cap I$
 and (if $n\ge 3$) $I_k=(x_{k-1},x_k)$, $k=2,\dots,n-1$.
 Let $x\in I_k$ for some $k\in\{1,\dots,n-1\}$. Fix $y\in I$ such that
 $y>x_{n-1}$. By Theorem~\ref{Th2} we get
 $[x_1,\dots,x_{n-1},x;f]_{\pmb{\omega}}\le[x_1,\dots,x_{n-1},y;f]_{\pmb{\omega}}$.
 Tending with $y$ to $x_{n-1}^+$ and using~\eqref{Th3_1} we obtain
 $[x_1,\dots,x_{n-1},x;f]_{\pmb{\omega}}\le c_n$, whence
 by~\eqref{div_diff}
 \begin{equation}\label{Th3_2}
  \frac{D_{n-1}(x_1,\dots,x_{n-1},x;f)}{V_n(x_1,\dots,x_{n-1},x)}\le c_n.
 \end{equation}
 We need $n-k$ inversions of~$x$ to transform the ordered system
 of~$n$ points $(x_1,\dots,x,\dots,$ $x_{n-1})$ to the system
 $(x_1,\dots,x_{n-1},x)$. Then
 \[
  0<V_n(x_1,\dots,x,\dots,x_{n-1})=(-1)^{n-k}V_n(x_1,\dots,x_{n-1},x).
 \]
 Hence multiplying both sides of an inequality~\eqref{Th3_2} by
 $(-1)^{n-k}V_n(x_1,\dots,$ $x_{n-1},x)$ we get
 \[
  (-1)^{n-k}\bigl(D_{n-1}(x_1,\dots,x_{n-1},x;f)-c_nV_n(x_1,\dots,x_{n-1},x)\bigr)\le 0
 \]
 and dividing both sides of this inequality by $V_{n-1}(x_1,\dots,x_{n-1})>0$
 we obtain~\eqref{Th3_i} (for $k=1$) and~\eqref{Th3_ii} (for
 $k=2,\dots,n-1$ if $n\ge 3$).
 \par
 Now we prove the sufficiency. Fix $x_1,\dots,x_{n+1}\in I$ such that
 $x_1<x_2<\dots<x_n<x_{n+1}$. By Theorem~A it is enough to check that
 $D_n(x_1,\dots,x_{n+1};f)\ge 0$. By the assumption there exist
 the constants $c_1,\dots,c_n\in\R$ such that for
 $\omega=c_1\omega_1+\dots+c_n\omega_n$ we have
 $\omega(x_k)=f(x_k)$, $k=2,\dots,n$
 and
 \begin{align}
  \label{Th3_4}
  f(x_{n+1})-\omega(x_{n+1})&\ge 0,\\
  \label{Th3_5}
  (-1)^n\bigl(f(x_1)-\omega(x_1)\bigr)&\ge 0.
 \end{align}
 Finally we expand the determinant $D_n(x_1,\dots,x_{n+1};f-\omega)$
 by the last row. By the definition of~$\omega$ its elements
 $f(x_k)-\omega(x_k)$ ($k=2,\dots,n$) are equal to zero. Since
 $\pmb{\omega}$ is a~positive Chebyshev system, the determinants
 $V_n(x_2,\dots,x_{n+1})$, $V_n(x_1,\dots,x_n)$ are positive.
 Since $D_n(x_1,\dots,x_{n+1};\omega)=0$ then
 by~\eqref{Th3_4},~\eqref{Th3_5} we infer
 \begin{multline*}
  D_n(x_1,\dots,x_{n+1};f)=
  D_n(x_1,\dots,x_{n+1};f-\omega)\\
  =(-1)^{n+2}\bigl(f(x_1)-\omega(x_1)\bigr)V_n(x_2,\dots,x_{n+1})\\
  +\bigl(f(x_{n+1})-\omega(x_{n+1})\bigr)V_n(x_1,\dots,x_n)\ge 0,
 \end{multline*}
 which finishes the proof.
\end{proof}
Using Theorem~\ref{Th3} for $n=2$ we obtain immediately the following
result (see~\cite[Theorem~4 $(i)\Leftrightarrow(iii)$]{BesPal03}).
\begin{cor}\label{Co2}
 Let $\pmb{\omega}=(\omega_1,\omega_2)$ be a~positive Chebyshev system
 on~$I$ such that $\omega_1>0$. A~function $f:I\to\R$ is
 $\pmb{\omega}$-$2$-convex if and only if for any $x_1\in I^0$ there
 exist the constants $c_1,c_2\in\R$ such that for
 $\omega=c_1\omega_1+c_2\omega_2$ we have $\omega(x_1)=f(x_1)$ and
 $\omega\le f$ on~$I$.
\end{cor}

\begin{rem}
 By Corollary~\ref{Co2} Theorem~\ref{Th3} reduces for $n=2$ to the
 classical support theorem. For $n\ge 3$ it is not the case. The
 function~$\omega$ supports~$f$ in the interval
 $(x_{n-2},+\infty)\cap I$. Passing through the points
 $\bigl(x_i,f(x_i)\bigr)$, $i=1,\dots,n-2$ the graph of~$\omega$
 successively changes the side of the graph of~$f$. Let us illustrate
 this situation by the following example.
\end{rem}

\begin{exmp}
 Let $n=3$ and $\pmb{\omega}=(1,x,x^2)$. Obviously $\pmb{\omega}$
 and $(1,x)$ are positive Chebyshev systems on any interval.
 By Theorem~A it is easy to see that $f(x)=x^3$ is
 $\pmb{\omega}$-$3$-convex ($D_3(x_1,x_2,x_3,x_4;f)$ is the
 Vandermonde determinant). Observe that the function
 $\omega(x)=2x^2-x$ fulfils the inequalities~\eqref{Th3_i},
 \eqref{Th3_ii} and~\eqref{Th3_iii} of Theorem~\ref{Th3} for
 $x_1=0$, $x_2=1$. Namely, $\omega(0)=f(0)$, $\omega(1)=f(1)$ and
 \begin{align*}
  f(x)-\omega(x)&\le 0\quad\text{for $x<0$},\\
  f(x)-\omega(x)&\ge 0\quad\text{for $0<x<1$},\\
  f(x)-\omega(x)&\ge 0\quad\text{for $x>1$}.
 \end{align*}
\end{exmp}
\begin{ack*}
 The author gratefully acknowledges the referee's remarks simplifying
 the proofs.
\end{ack*}

\bibliographystyle{amsplain}
\bibliography{was_pub}
\end{document}